\documentclass[final,3p,times]{elsarticle}
\usepackage{amssymb}
\usepackage{graphicx,setspace}
\usepackage{psfrag}
\usepackage{amsfonts}
\expandafter\let\csname equation*\endcsname\relax
\expandafter\let\csname endequation*\endcsname\relax
\usepackage{amsmath}
\usepackage{amssymb, amsthm}
\usepackage{mathrsfs}
\usepackage{epstopdf}
\usepackage{float}
\usepackage{lineno,hyperref}
\usepackage{color}
\usepackage{subfigure}
\usepackage{amsmath}
\usepackage{dsfont} 
\usepackage{amssymb} 
\usepackage{graphicx,color}
\usepackage{extarrows}

\usepackage{graphicx}

\usepackage{epstopdf}

\journal{arXiv}

\begin{document}
\newtheorem{definition}{Definition}[section]
\newtheorem{lemma}{Lemma}[section]
\newtheorem{remark}{Remark}[section]
\newtheorem{theorem}{Theorem}[section]
\newtheorem{proposition}{Proposition}
\newtheorem{assumption}{Assumption}
\newtheorem{example}{Example}
\newtheorem{corollary}{Corollary}[section]
\def\ep{\varepsilon}
\def\Rn{\mathbb{R}^{n}}
\def\Rm{\mathbb{R}^{m}}
\def\E{\mathbb{E}}
\def\hte{\hat\theta}
\renewcommand{\theequation}{\thesection.\arabic{equation}}
\begin{frontmatter}

\title{Endogenous business cycles via state-dependent saving and noise-induced metastability }

\author{Shenglan Yuan\corref{cor1}\fnref{addr1}}\ead{shenglanyuan@gbu.edu.cn}\cortext[cor1]{Corresponding author}

\author{James Brannan\fnref{addr2}}\ead{jrbrn@clemson.edu}

\author{Almaz Abebe\fnref{addr3}}\ead{almaz.abebe@howard.edu}

\author{Daniel Tesfay\fnref{addr4}}\ead{daniel.tesfay@chalmers.se}

\address[addr1]{\rm Department of Mathematics, School of Sciences, Great Bay University, Dongguan 523000, China }
\address[addr2]{\rm Department of Mathematical and Statistical Sciences, Clemson University, Clemson, South Carolina 29634, USA}
\address[addr3]{\rm Department of Mathematics, College of Art and Science, Howard University, Washington DC, 20059, USA}
\address[addr4]{\rm Department of Mathematical Sciences,
Chalmers University of Technology,
41296 Gothenburg, Sweden}

\begin{abstract}
We develop a parsimonious stochastic growth model in which
state-dependent saving behavior generates endogenous
business-cycle-like dynamics.
The model consists of three coupled equations:
a Solow-type capital accumulation equation,
a linear filtering equation for the saving rate,
and a bounded stochastic adjustment process.
Saving is modeled as a logistic function of deviations
from a balanced growth path,
introducing nonlinear feedback controlled by a gain parameter.

In the deterministic limit, increasing feedback strength
produces a supercritical pitchfork bifurcation,
splitting the balanced-growth equilibrium into two locally
attracting regimes corresponding to expansion and contraction.
When stochastic perturbations are introduced,
these equilibria become metastable states,
and the economy undergoes rare noise-induced transitions
between them.
The resulting dynamics exhibit persistent regimes,
bimodal stationary densities, and right-skewed dwell-time
distributions with approximately exponential survival tails.

A discrete-time approximation is estimated using U.S. real GDP data,
and Monte Carlo simulations are used to compute stationary
distributions and regime persistence statistics.
The results demonstrate that nonlinear state dependence,
bounded multiplicative noise, and time-scale separation
are sufficient to generate realistic business-cycle behavior
within a low-dimensional framework.
\end{abstract}

\begin{keyword}
Endogenous business cycles, State-dependent saving, Stochastic bifurcation,  Metastability

\emph{2020 Mathematics Subject Classification}: 60H30, 37H20

\end{keyword}

\end{frontmatter}

\section{Introduction}
Real GDP exhibits persistent long-run growth together with
recurrent short-run fluctuations commonly referred to as business cycles \cite{CG,CFS,FB,KG,KM}.
While a large literature explains these fluctuations through
exogenous shocks, frictions, or highly structured multi-sector models \cite{A,B},
a fundamental question remains:
to what extent can business-cycle-like dynamics arise endogenously
from simple nonlinear feedback mechanisms operating under uncertainty?

This work develops a minimal stochastic growth model in which
state-dependent saving behavior, coupled with bounded stochastic
perturbations, generates regime persistence and random transitions
between expansionary and contractionary states.
The emphasis is not on reproducing every feature of macroeconomic data,
but on identifying a transparent geometric mechanism capable of producing
business-cycle phenomena within a low-dimensional system.

The starting point is a Solow-type capital accumulation equation \cite{S}
formulated in terms of a dimensionless residual variable $k(t)$,
which measures deviations of capital from a locally defined
balanced growth path.
The key modeling ingredient is a nonlinear saving function
that depends on the state of the economy.
Saving is interpreted broadly as the aggregate share of output
devoted to capital formation, including corporate investment,
government capital expenditures, and household saving
intermediated into productive capital.

The saving rate is modeled as a bounded logistic function of $k$,
with a gain parameter controlling the sensitivity of saving
to deviations from trend.
This nonlinear feedback introduces the possibility of
multiple deterministic equilibria.
As the gain parameter increases, the system undergoes
a supercritical pitchfork bifurcation in which the balanced-growth
equilibrium loses stability and two locally attracting regimes emerge.
This deterministic bifurcation provides the geometric backbone
of the model.

When stochastic fluctuations are introduced into the saving dynamics,
the two deterministic attractors become metastable states.
The system spends extended periods fluctuating near one regime
before undergoing rare noise-induced transitions to the other.
The resulting dynamics resemble business cycles:
persistent expansions and contractions with random durations.
In contrast to models driven by periodic forcing,
the cycles here arise from stochastic transitions between
nonlinearly generated regimes.

The model possesses a natural fast--slow structure.
A rapidly adjusting stochastic component governs short-run
fluctuations in saving behavior, while the capital stock evolves
on a slower time scale.
A linear filtering equation connects these two components,
producing smooth saving dynamics despite high-frequency noise.
This time-scale separation plays a central role in shaping
both the deterministic bifurcation and the stochastic regime-switching behavior.

To connect the model with data,
we construct a discrete-time approximation suitable for estimation
using U.S. real GDP.
Parameters are estimated via a prediction-error framework
augmented by moment-matching constraints.
Monte Carlo simulations are then used to compute stationary
distributions, transition statistics, and dwell-time distributions.
The simulated dwell times exhibit right-skewed distributions
with approximately exponential survival tails,
consistent with the theory of noise-induced escape from metastable states.

The contributions of the work are threefold.
First, we introduce a parsimonious stochastic growth model
in which state-dependent saving generates endogenous regime switching
through a clearly identifiable bifurcation mechanism.
Second, we show how bounded multiplicative noise interacting
with nonlinear feedback produces metastability and realistic
cycle-length statistics without requiring large exogenous shocks.
Third, we provide a simulation-based framework linking
geometric dynamical systems theory to empirical macroeconomic data.

The remainder of the paper is organized as follows.
Section~2 motivates the state-dependent saving specification
and the time-scale separation underlying the model.
Section~3 presents the continuous-time system and its
fast--slow geometric structure.
Section~4 analyzes the deterministic bifurcation
and develops a discrete-time approximation for estimation.
Section~5 reports Monte Carlo results on stationary densities
and regime persistence.
Section~6 concludes.
\section{State-Dependent Saving and Time-Scale Separation}

The central modeling assumption of this work is that the aggregate saving
rate depends nonlinearly on the state of the economy and evolves on a
faster time scale than the capital stock.
This section provides the economic and empirical motivation for these
two assumptions.

\subsection{Logistic Feedback in the Saving Rate}

We model the state-dependent saving rate by the logistic function
\begin{equation}
g(k;\gamma)
=
s_1
+
\frac{s_2 - s_1}{1 + e^{-\gamma (k-1)}},
\label{logistic_s}
\end{equation}
where $0 \le s_1 < s_2 \le 1$ and $\gamma > 0$ controls the sensitivity
of saving behavior to deviations of capital from its reference level $k=1$.

The logistic form captures three economically plausible regimes:

\begin{itemize}
\item[(i)] \textbf{Low capital ($k \ll 1$).}
At low levels of capital, income is close to subsistence.
A large fraction of output is devoted to consumption,
and the saving rate approaches its lower bound $s_1$.

\item[(ii)] \textbf{Intermediate capital.}
As capital and productivity rise,
households and firms generate surplus resources.
The saving rate increases more rapidly,
with maximal sensitivity near $k=1$.

\item[(iii)] \textbf{High capital ($k \gg 1$).}
At sufficiently high income levels,
behavioral, institutional, or policy constraints
limit further increases in the investment share.
The saving rate saturates near $s_2$.
\end{itemize}

The logistic specification provides a smooth transition between
low- and high-saving regimes while maintaining boundedness.
More importantly, its slope near equilibrium is controlled by $\gamma$.
As $\gamma$ increases, the feedback from capital to saving becomes sharper,
introducing the possibility of nonlinear amplification and,
as shown below, bifurcation in the deterministic dynamics.

\subsection{Time-Scale Separation and Stochastic Fluctuations}

The second modeling ingredient is that saving behavior adjusts more
rapidly than the capital stock.
Capital accumulation reflects investment decisions, production capacity,
and structural adjustment, all of which evolve gradually.
In contrast, saving behavior responds quickly to expectations,
financial conditions, policy changes, and short-run sentiment.

Empirically, the personal saving rate is observed at monthly frequency
and exhibits substantial short-run variability,
whereas real GDP evolves more smoothly and is measured quarterly.
This contrast supports modeling saving as a higher-frequency variable
relative to capital.

Accordingly, we introduce an intermediate variable $p(t)$ that
reverts rapidly toward the state-dependent target $\eta g(k(t))$,
while remaining confined to a bounded interval.
The observable saving rate $s(t)$ then evolves according to a stable
relaxation equation driven by $p(t)$.
This structure produces three distinct time scales:

\begin{itemize}
\item a \emph{fast} stochastic adjustment of $p(t)$ toward $g(k)$,
\item an \emph{intermediate} filtering dynamics for $s(t)$,
\item a \emph{slow} evolution of the capital stock $k(t)$.
\end{itemize}

The resulting system combines nonlinear state dependence,
bounded stochastic perturbations, and time-scale separation.
As shown in the following sections, this minimal structure is
sufficient to generate multiple deterministic equilibria and,
in the presence of noise, stochastic regime switching.

\section{The Mathematical Model}
Consider the following system,
\begin{align}
d\tilde{k}(t) &= \big[\tilde{s}(t)\,b\,\tilde{k}(t)^{\alpha} - \mu\,\tilde{k}(t)\big]\, dt,
\label{dimensioned_k_equation}\\
d\tilde{s}(t) &= -\beta \big[\tilde{s}(t)- \tilde{p}(t)\big]\, dt, \label{dimensioned_s_equation} \\
d\tilde{p}(t) &= -a\big[\tilde{p}(t) - \tilde{g}(\tilde{k}(t))\big]\, dt
\;+\; \sigma\,\sqrt{\tilde{p}(t)\big(1 - \tilde{p}(t)\big)}\, dW(t), \label{dimensioned_p_equation}
\end{align}
where
\[
\tilde{g}(\tilde{k}) = s_1 + \frac{s_2 - s_1}{1 + e^{-\tilde{\gamma}(\tilde{k}-\tilde{k}^{\ast})}},
\qquad
\tilde{k}^{\ast} = \left( \frac{b(s_1 + s_2)}{2 \mu} \right)^{1/(1-\alpha)} .
\]

\subsection{Scaled Variables}

Introducing the dimensionless state variables $k(t)$, $s(t)$, and $p(t)$
using the scalings
\[
\tilde{k}(t) = \tilde{k}^{\ast}\,k(t),
\qquad
\tilde{s}(t) = \frac{s_1+s_2}{2}\,s(t),
\qquad
\tilde{p}(t) = \frac{s_1+s_2}{2}\,p(t)
\]
yields the following system
\begin{align}
dk(t) &= \mu\big[ s(t)\, k(t)^{\alpha} - k(t) \big]\, dt,
\label{scaled_k_equation} \\
ds(t) &= -\beta \big[s(t)-p(t)\big] dt, \label{scaled_s_equation} \\
dp(t) &= -a \big[p(t)- \eta g(k(t))\big] dt + \sigma\sqrt{p(t)(\eta-p(t))}\,dW(t),
\label{scaled_p_equation}
\end{align}
where
\[
\eta = \frac{2}{s_1+s_2}, \qquad
\gamma = \tilde{\gamma}\,\tilde{k}^{\ast}, \qquad
\]
and
\begin{equation}\label{SigmoidFunction}
g(k(t)) = s_1 + \frac{s_2-s_1}{1+ e^{-\gamma(k(t)-1)}}.
\end{equation}
Note that the parameters $\mu$,  $\beta$, and $a$
have dimensions of $time^{-1}$ while the dimensions of $\sigma$ are
$time^{-1/2}$.

The logistic function $g(k)$ describes an endogenous savings rate
bounded between $s_1$ and $s_2$, with responsiveness controlled by
$\gamma$.  The centering point $k^*$ is chosen as the positive
equilibrium of the model with constant savings rate
$\bar{s}=(s_1+s_2)/2$:
\begin{equation}
  \bar{s}\,b\,(k^*)^\alpha = \mu k^*,
  \qquad
  k^* =
  \left(\frac{b(s_1+s_2)}{2\mu}\right)^{\!\frac{1}{1-\alpha}}.
\end{equation}

\subsection{Geometric Explanation of the Model}

Equation (\ref{scaled_k_equation}) is an intensive Solow-type equation in which
$k(t)$ denotes the residual capital stock per unit of effective labor,
measured relative to a nearby balanced growth path.
The saving rate $s(t)$ represents the aggregate fraction of output
devoted to capital formation.
It includes corporate capital expenditures, government capital expenditures,
and household saving that is intermediated into productive investment.
Thus $s(t)Y(t)$ equals total gross investment in the economy.

The full system (\ref{scaled_k_equation})--(\ref{scaled_p_equation}) possesses a natural fast--slow structure.
For fixed $k$, the variable $p(t)$ evolves according to the
Wright--Fisher-type diffusion (\ref{scaled_p_equation}), with drift
$-a\big(p-\eta g(k)\big)$ and diffusion coefficient
$\sigma \sqrt{p(\eta-p)}$.
The drift term pulls $p(t)$ toward the logistic curve
\[
p=\eta g(k),
\]
while the multiplicative noise keeps the process confined to $(0,\eta)$.
When $a \gg 1$, the relaxation of $p(t)$ toward this curve occurs
on a time scale much shorter than that of $s(t)$ and $k(t)$.

In the singular limit $a \to \infty$, Eq. (\ref{scaled_p_equation}) reduces to the
algebraic constraint
\[
p=\eta g(k),
\]
and trajectories collapse onto the one-dimensional critical manifold
\[
\mathcal{M}_0=\{(k,s,p)\in\mathbb{R}^3 : p=\eta g(k)\}.
\]
For large but finite $a$, geometric singular perturbation theory
implies the existence of a locally invariant slow manifold
\[
\mathcal{M}_a=\{(k,s,p): p=\eta g(k)+\mathcal{O}(a^{-1})\},
\]
which attracts nearby trajectories exponentially fast.
After a short transient, the dynamics evolve within a thin stochastic
tube surrounding this manifold in $(k,s,p)$-space.

The saving rate $s(t)$ does not instantaneously equal $p(t)$.
Instead, it satisfies the linear relaxation equation
\[
\frac{ds}{dt} = -\beta (s-p),
\]
whose unique bounded solution is
\[
s(t)=\beta \int_0^\infty e^{-\beta \zeta}\, p(t-\zeta)\, d\zeta.
\]
This representation shows that $s(t)$ is an exponentially weighted
moving average of $p(t)$.
In signal-processing terminology, equation (\ref{scaled_s_equation}) acts as a
stable low-pass filter:
high-frequency fluctuations of $p(t)$ are attenuated, while
low-frequency components are transmitted.

Geometrically, $p(t)$ undergoes rapid stochastic motion within a
thin neighborhood of the slow manifold $p=\eta g(k)$, whereas
$s(t)$ evolves as a smoothed projection of this motion onto
a slower time scale.
Consequently, although the fast variable $p(t)$ exhibits
Wright--Fisher-type fluctuations, the observable saving rate $s(t)$
varies more smoothly and remains confined to $(0,\eta)$.
It is this filtered feedback of $s(t)$ into the capital equation
(\ref{scaled_k_equation}) that produces nonlinear amplification and, in the
deterministic limit, the pitchfork bifurcation described below.

The model is deliberately minimal: three coupled equations with
bounded stochastic forcing suffice to generate regime persistence
and stochastic transitions \cite{ATA,YH} between expansionary and contractionary
states.
The richness of the resulting dynamics arises not from high dimension,
but from the interaction of nonlinear state dependence,
time-scale separation, and noise confined to a compact interval.

\subsection{Deterministic Dynamics and Bifurcation Structure}

Setting $\sigma = 0$ in Eqs. (\ref{scaled_k_equation})--(\ref{scaled_p_equation}) yields the deterministic system
\begin{equation}
\frac{dy}{dt} = H(y;\gamma),
\qquad
y =
\begin{bmatrix}
k \\ s \\ p
\end{bmatrix},
\label{det_system}
\end{equation}
where
\[
H(y;\gamma) =
\begin{bmatrix}
\mu(s k^\alpha - k) \\
-\beta (s - p) \\
-a\big(p - \eta g(k;\gamma)\big)
\end{bmatrix}.
\]

Equilibrium solutions satisfy
\[
p = \eta g(k;\gamma),
\qquad
s = p,
\qquad
s k^\alpha - k = 0.
\]
Eliminating $s$ and $p$ yields the scalar equation
\begin{equation}
\eta g(k;\gamma) - k = 0.
\label{scalar_equilibrium}
\end{equation}

The centering choice of scaling ensures that $k=1$ solves
\eqref{scalar_equilibrium} for all $\gamma \ge 0$, giving the equilibrium
\[
y_0^*(\gamma) = (1,1,1)^T.
\]

\paragraph{Linearization and Loss of Stability}

Linearizing \eqref{det_system} about $y_0^*(\gamma)$ gives
\[
\frac{du}{dt} = A(\gamma)u + \mathcal{O}(\|u\|^2),
\]
where
\begin{equation}\label{A_MatrixFromLinearization}
A(\gamma) =
\begin{bmatrix}
\mu(\alpha - 1) & \mu & 0 \\
0 & -\beta & \beta \\
\frac{a\eta\gamma (s_2 - s_1)}{4} & 0 & -a
\end{bmatrix}.
\end{equation}

The characteristic polynomial of $A(\gamma)$ has constant term
\[
\det A(\gamma)
=
\mu a \beta
\left[
(\alpha - 1) + \frac{\eta \gamma (s_2 - s_1)}{4}
\right].
\]

A change in stability occurs when $\det A(\gamma)=0$, yielding the critical value
\begin{equation}
\gamma_c
=
\frac{2(1-\alpha)}{s_2+s_1}
\frac{s_2+s_1}{s_2 - s_1}
=
\frac{2(1-\alpha)}{s_2 - s_1}(s_1+s_2).
\label{gamma_c}
\end{equation}

At $\gamma=\gamma_c$, a simple real eigenvalue crosses zero transversely,
while the remaining eigenvalues remain strictly negative.
Thus the equilibrium $y_0^*(\gamma)$ loses stability through a
codimension-one bifurcation.

\paragraph{Pitchfork Structure}

To determine the nature of the bifurcation, it is sufficient to examine
the scalar reduced equation \eqref{scalar_equilibrium}.
The function
\[
F(k;\gamma) := \eta g(k;\gamma) - k
\]
is odd about $k=1$ due to the symmetric centering of the logistic function,
and satisfies
\[
\partial_k F(1;\gamma_c)=0,
\qquad
\partial_{kk} F(1;\gamma_c)=0,
\qquad
\partial_{kkk} F(1;\gamma_c)<0.
\]

Hence, at $\gamma=\gamma_c$, the equilibrium at $k=1$ undergoes a
supercritical pitchfork bifurcation.
For $\gamma < \gamma_c$, $k=1$ is the unique equilibrium and is asymptotically stable.
For $\gamma > \gamma_c$, two additional equilibria
$k_\pm(\gamma)=1 \pm \delta(\gamma)$ emerge symmetrically,
with $k_+(\gamma)>1$ and $k_-(\gamma)<1$, both stable,
while $k=1$ becomes unstable.

\paragraph{Geometric Interpretation}

The bifurcation is driven entirely by the gain parameter $\gamma$,
which controls the sensitivity of saving behavior to deviations of
capital from trend.
When $\gamma$ is small, the feedback from $k$ to $s$ is weak,
and the restoring force toward the balanced growth path dominates.
As $\gamma$ increases, the slope of the sigmoid near $k=1$ steepens,
amplifying deviations in saving and thereby strengthening positive feedback.

\begin{figure}[H]
\begin{center}
\includegraphics[width=12cm]{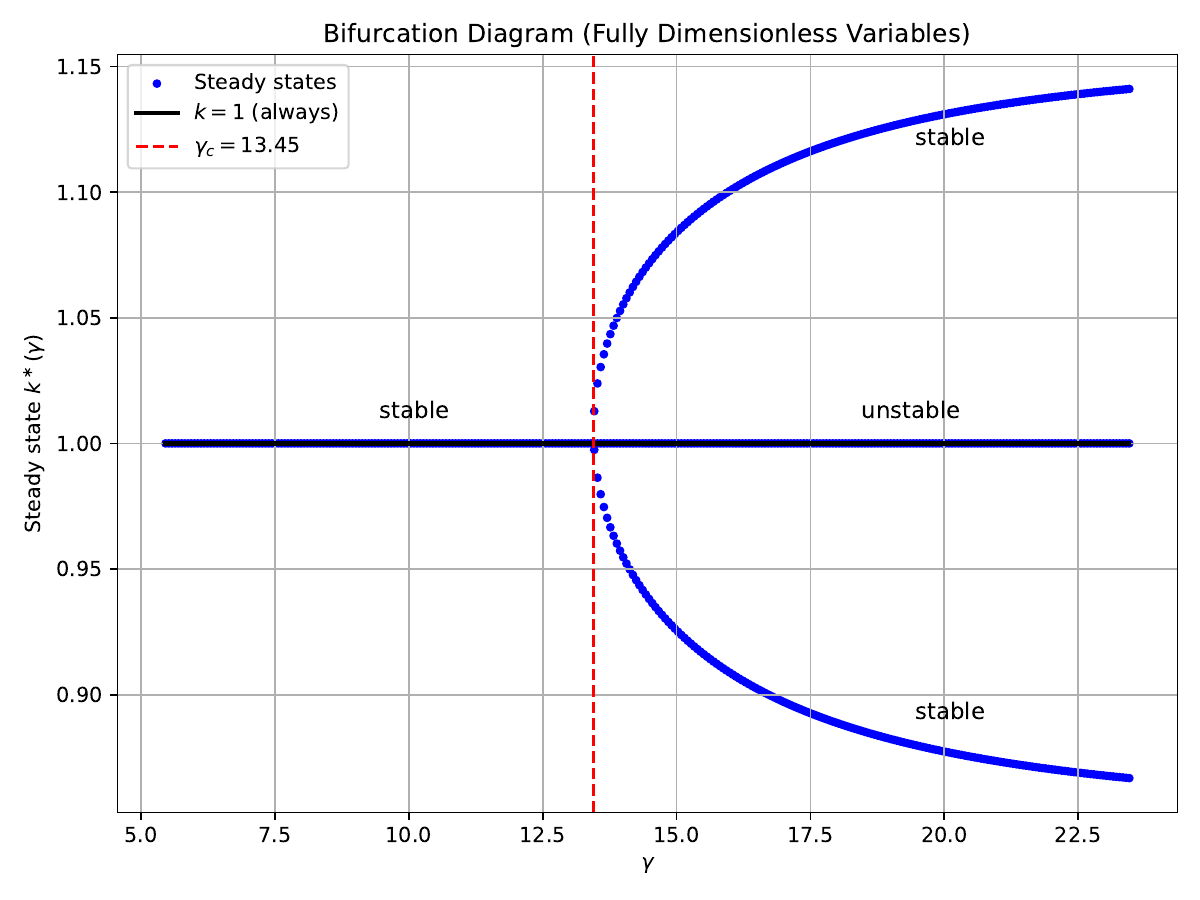}
\caption{Bifurcation diagram for the first component of Eq. (\ref{det_system})}.
\label{BifurcationDiagram1}
\end{center}
\end{figure}

At the critical value $\gamma_c$, the stabilizing linear term in the
capital equation is exactly offset by the nonlinear feedback through
the saving function.
Beyond this threshold, the balanced growth equilibrium splits into
two locally attracting regimes corresponding to persistent expansion
($k>1$) and contraction ($k<1$).
The bifurcation diagram for the first component of equilibrium solutions of
Eq. (\ref{det_system}) is shown in Figure \ref{BifurcationDiagram1}.
This deterministic pitchfork structure forms the geometric backbone
of the stochastic model.
In the presence of noise, the two stable deterministic equilibria
become metastable states between which the system transitions randomly.

\section{A Discrete Approximation and Its Identification}\label{DiscreteSection}
We obtain a discrete approximation to system (\ref{scaled_k_equation})-(\ref{scaled_p_equation}) by setting $k_n = k(n\Delta t)$, $s_n = s(n\Delta t)$,
$p_n = p(n\Delta t)$, replacing time derivatives by difference quotients, and absorbing $\Delta t$ into the parameters:
\vspace{20pt}
\begin{align}
k_{n+1} & =  k_n + \mu \big(s_{n+1} k_n^\alpha-k_n \big),  \label{discrete_k_eqn} \\
s_{n+1} & =  s_n - \beta \big(s_n-p_{n+1} \big), \label{discrete_s_eqn}\\
p_{n+1} & = p_n - a\big(p_n-\eta g(k_n)\big)
+ \sigma \sqrt{p_n(\eta-p_n)}\, \varepsilon_{n+1}, \label{discrete_p_eqn}
\end{align}
where
\[
g(k_n)
= s_1 + \frac{s_2 - s_1}{1 + e^{-\gamma (k_n - 1)}}
\qquad \text{and} \qquad
\{\varepsilon_n\} \stackrel{\text{i.i.d.}}{\sim} \mathcal{N}(0,1).
\]
Discrete time $n$ in system (\ref{discrete_k_eqn})-(\ref{discrete_p_eqn}) refers to quarters of a year (periods of three months) in real time.  Thus $n=4$ corresponds to 1 year of elapsed time, $n=6$ corresponds to $1\frac{1}{2}$
years of elapsed time, and so on.  Let us denote the parameters in the above system by the vector
$\theta~=~\left(\alpha, s_1, s_2, \mu, \gamma, \beta, a, \sigma \right)$.
Note that $\eta=(2/(s_1+s_2)$ is not an independent parameter.
Occasionally it will be convenient
to use the vector notation $\mathbf{y}_n = (k_n, s_n, p_n)^T$ to represent the state of the system.
\subsection{Bifurcation Analysis of the Discrete Deterministic Model}\label{DiscreteBifurcationAnalysis}
\begin{figure}[H]
\begin{center}
\includegraphics[width=7.7cm]{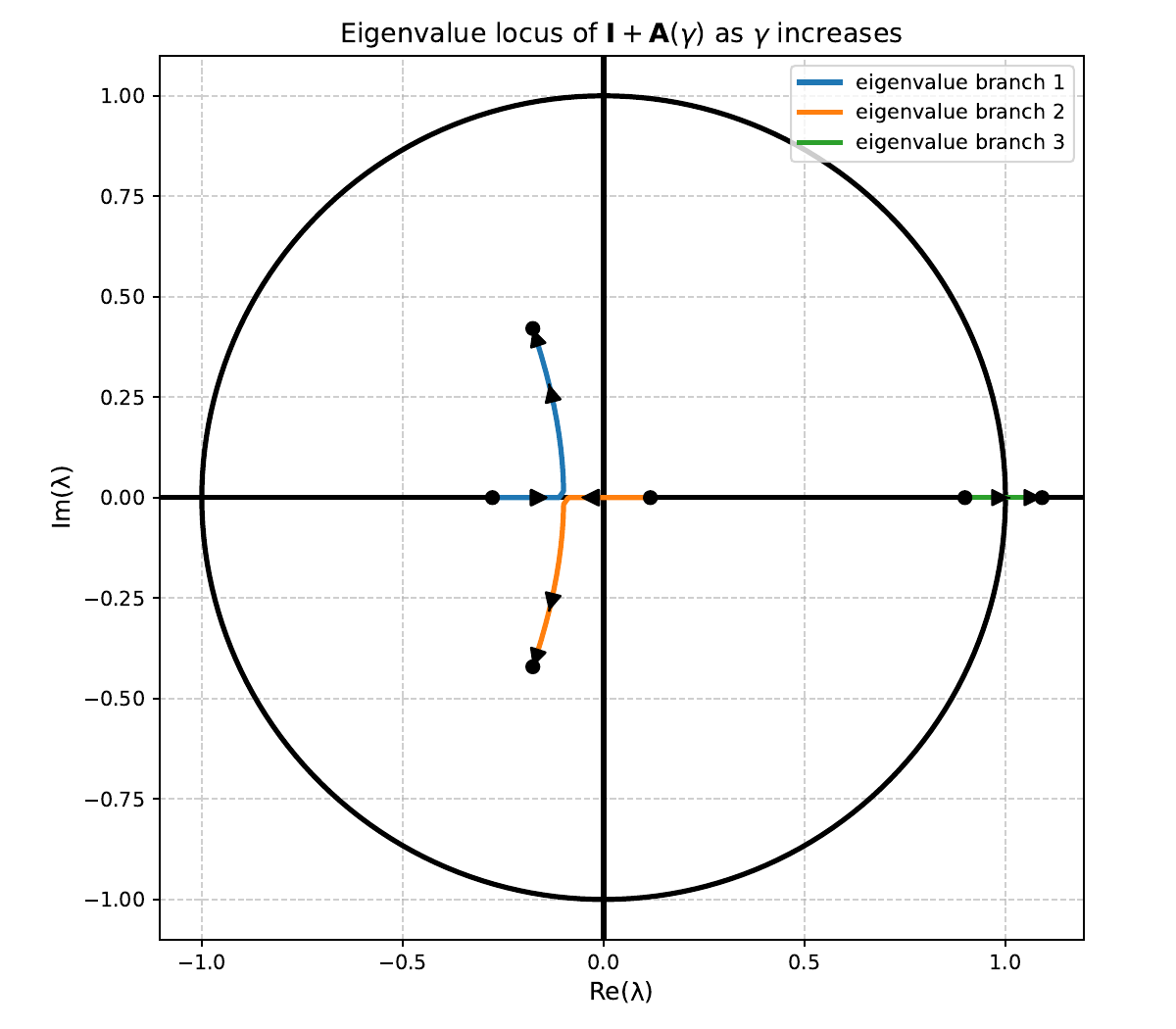}
\caption{Loci of the eigenvalues of $\mathbf{I}+\mathbf{A}(\gamma)$ in
Eq. (\ref{DiscreteLinearSystem}) as $\gamma$ increases from values
below $\gamma_c$ to values above $\gamma_c$.}
\label{BifurcationDiagram}
\end{center}
\end{figure}

Before we estimate the parameters in the system (\ref{discrete_k_eqn})-(\ref{discrete_p_eqn}) using GDP data in Section \ref{PredictionErrorEstimation}, it is of interest to examine the steady state solutions of the system.  Linearizing
Eqs. (\ref{discrete_k_eqn})-(\ref{discrete_p_eqn}) around the equilibrium solution
\begin{equation}
\mathbf{x}_0^{\ast} =
\left[
\begin{array}{c}
  1   \\
  1   \\
  1
\end{array}
\right]
\end{equation}
 gives
 \begin{equation}\label{DiscreteLinearSystem}
\mathbf{u}_{n+1} = \Big(\mathbf{I} + \mathbf{A}(\gamma)\Big)  \mathbf{u}_{n},
\end{equation}
where $\mathbf{A}(\gamma)$ is the matrix in Eq. (\ref{A_MatrixFromLinearization}).  If $\lambda = 0$ is an eigenvalue of
$\mathbf{A}(\gamma_c)$, $\lambda = 1$ is an eigenvalue of
$\mathbf{I} + \mathbf{A}(\gamma_c)$.  As shown in Figure \ref{BifurcationDiagram}, one of the
eigenvalues of $\mathbf{I}+\mathbf{A}(\gamma)$ moves along the
real axis from inside the unit circle to the exterior of the unit circle, taking on
the value 1 when $\gamma=\gamma_c$.

\subsection{Prediction Error Model Estimation}\label{PredictionErrorEstimation}
Denote by $\left\{ (t_n, Y(t_n)), \; n = 1,2, \ldots, N \right\}$ real GDP data where $t_n = 0.25n$ years.\footnote{This data is posted by the U.S. Bureau of Economic Analysis, Real Gross Domestic Product [GDPC1], retrieved from FRED, Federal Reserve Bank of St. Louis;
https://fred.stlouisfed.org/series/GDPC1, November 18, 2026.  The $t_n$ are the quarterly dates of record while the $Y(t_n)$ are the corresponding seasonally adjusted values of the GDP in billions of chained 2012 dollars.}
We derive the sequence of data points
$\displaystyle
\left\{ (\hat{k}_n, \hat{s}_n, \hat{p}_n, \hat{\varepsilon}_n), \;n=1,2, \ldots,N \right\}
$
from the time series $\left\{ Y(t_n), \; n = 1,2, \ldots, N \right\}$ as follows.
From the assumed model for the GDP
\[
Y(t_n) = A_0L_0e^{(n+g)t_n}Bk^{\alpha}(t_n),
\]
we get
\[
\hat{k}_n = k(t_n) = \exp\left\{ \frac{1}{\alpha} \left[ \ln y_n - \ln{(A_0L_0B\phi^{\alpha})} - (n+g) t_n \right] \right\},
\]
where
\[
y_n = Y(t_n), \quad n = 1,2, \ldots, N.
\]
Using $t_n=0.25n$ years as the independent variable, linear regression on the time series  $\left\{\ln(y_n), n = 1, 2, \ldots, N \right\}$ gives us
\begin{equation}\label{kndata2}
\hat{k}_n = \hat{k}(t_n) = \exp\left\{ \frac{1}{\alpha} \left[ \ln y_n - b_0 - b_1 t_n \right] \right\},
\end{equation}
where $b_0$ is the estimate of $\ln{(A_0L_0B\phi^{\alpha})}$ and $b_1$ is the estimate of $(n+g)$.

Given the data $\left\{ \hat{k}_n, \; n = 1, \ldots N \right\}$, assuming that $s_1=1$, the rest of the sequence $\left\{ \hat{s}_n, \; n = 2, \ldots N \right\} $ is obtained from Eq. (\ref{discrete_k_eqn}),
\[
\hat{s}_{n+1} = \frac{\hat{k}_{n+1}-(1-\mu)\hat{k}_n}{\mu \hat{k}_n^{\alpha}},
\qquad n = 2, \ldots,N-1.
\]
Next, assume that $\hat{p}_1 = 1$ and then solve Eq. (\ref{discrete_s_eqn}) for $p_n$ to generate the sequence
\[
\hat{p}_{n+1} = \frac{\hat{s}_{n+1} - (1-\beta) \hat{s}_n}{\beta}, \quad n =  2, \ldots, N-1.
\]
Finally, assume that $\hat{\varepsilon}_1 = 0$ and obtain the sequence $\left\{ \hat{\varepsilon}_n, \; n = 1, \ldots N \right\}$ by solving
Eq. (\ref{discrete_p_eqn}) for $\varepsilon_n$,
\begin{equation}\label{x_ndata}
\hat{\varepsilon}_{n+1} = \frac{\hat{p}_{n+1}-\hat{p}_n+
a[\hat{p}_n-\eta g(\hat{k}_n)]}{\sigma \sqrt{\hat{p}_n(\eta-\hat{p}_n)}}
 \quad n = 1, \ldots, N-1.
\end{equation}
Equation (\ref{x_ndata}) is where fluctuations in the GDP about the balanced
growth path are transformed into fluctuations represented by the time series
$\left\{\hat{\varepsilon}_{n}, \; n= 1, 2,\ldots,  N \right\}$.

At this point we review a few theoretical constructs from probability theory to explain the method we use to estimate the parameters in system
(\ref{discrete_k_eqn})-(\ref{discrete_p_eqn}).
Given a stochastic process, say $\left\{ k_n, \; n = 1, 2, \ldots \right\}$,
we denote by $\mathcal{F}_n=\langle k_j, \; 1 \le j \le n \rangle$, the sigma algebra generated by the random variables $k_1, k_2, \ldots,k_n$.
The collection $\mathcal{F}_n$  of sub $\sigma$-algebras is referred to as a \emph{filtration} since for every $m \leq n$, we have $\mathcal{F}_m \subseteq \mathcal{F}_n$.  The filtration $\mathcal{F}_n, n = 1, 2, \cdots$ models the information that is available at the given times $n=1, 2, \cdots$.  This definition generalizes to vector valued stochastic processes.  Treating $\left\{(k_n, s_{n}, p_{n},\varepsilon_{n}), \; n = 1,2, \cdots \right\}$ as a sequence of random vectors, we denote by $\mathcal{G}_n=\langle (k_j, s_{j}, p_{j},\varepsilon_{j}), \; 1\le j \le n \rangle$ the sigma algebra generated by $\left\{(k_j, s_{j}, p_{j},\varepsilon_{j}), \; 1 \le j \le n \right\}$.
 Obviously
\[
\mathcal{F}_{n} = \mathcal{G}_n,
\]
that is, the information in $\langle (k_j, z_j, p_{j},x_{j}), \; j \le n \rangle$
is identical to the information in $\langle k_n, \; j \le n \rangle$.
To obtain what is referred to as a one-step-ahead predictor of $k_{n+1}$ given $k_1, \ldots, k_n$, we find the function $H(k_1, \ldots, k_n)$ that minimizes
$E[\left\{ k_{n+1} - H(k_1, \ldots, k_n)\right\}^2]$.  From probability theory it is known that $H(k_1, \ldots, k_n)$ is given by the conditional expectation $E_{\theta}[k_{n+1}|\mathcal{F}_{n}]$, where $\theta$ indicates the dependence on the parameter vector
$\theta$.

Using the facts that $E_{\theta}[k_n|\mathcal{F}_{n}] = k_n$,
$E_{\theta}[s_{n}|\mathcal{F}_{n}] = s_{n}$, $E_{\theta}[p_{n}|\mathcal{F}_{n}] = p_{n}$, and $E_{\theta}[\varepsilon_{n+1}|\mathcal{F}_{n}] = 0$
we compute $E_{\theta}[k_{n+1}|\mathcal{F}_{n}]$ from the state space
model (\ref{discrete_k_eqn})-(\ref{discrete_p_eqn}) using the equations

\begin{gather*}
p_n
= p_{n-1} - a\bigl[p_{n-1} - \eta g(k_{n-1})\bigr]
+ \sigma \sqrt{p_{n-1}\bigl(\eta - p_{n-1}\bigr)}\,\varepsilon_n, \\
E_{\theta}\!\left[ p_{n+1} \mid \mathcal{F}_n \right]
= (1-a)p_n + a\eta g(k_n), \\
E_{\theta}\!\left[ s_{n+1} \mid \mathcal{F}_n \right]
= (1-\beta)s_n + \beta\,E_{\theta}\!\left[ p_{n+1} \mid \mathcal{F}_n \right], \\
E_{\theta}\!\left[ k_{n+1} \mid \mathcal{F}_n \right]
= k_n + \mu \left\{ E_{\theta}\!\left[ s_{n+1} \mid \mathcal{F}_n \right] k_n^{\alpha} - k_n \right\}.
\end{gather*}

The prediction error cost function is given by
\[
J_N^{(1)}(\theta) = \frac{1}{N-1} \sum_{n=2}^{N} \left\{E_{\theta}[k_{n}|\mathcal{F}_{n-1}]-\hat{k}_{n} \right\}^2,
\]
where the observed value of $k_n$ is denoted by $\hat{k}_n$.
In addition we impose constraints on the sample mean, sample standard deviation, and sample autocorrelation of a large sample simulation $\left\{k_n, \;n=1,\ldots, M,\; M \gg N \right\}$ of the system (\ref{discrete_k_eqn})-(\ref{discrete_p_eqn}) as follows.
Let
\[
g_1(\theta) = \frac{1}{M}\sum_{n=1}^M k_n,
\]
\[
g_2(\theta) = \left\{{\frac{1}{M-1}\sum_{n=1}^M \left[k_n-g_1(\theta]\right)^2} \right\}^{1/2},
\]
and
\[
g_3(m,\theta) =\frac{1}{M-m}\sum_{n=1}^M \left[k_{n+m}-g_1(\theta))(k_n-g_1(\theta)\right]/g_2^2(\theta)
\]
be estimates of the mean, standard deviation, and autocorrelation of the sequence $\left\{k_n, \;n=1,\ldots, M\right\}$ generated from the system (\ref{discrete_k_eqn})-(\ref{discrete_p_eqn}) and let
\[
\hat{\kappa} = \frac{1}{N} \sum_{n=1}^N \hat{k}_n, \quad
\hat{\chi} =\left\{  \frac{1}{N-1} \sum_{n=1}^N \left( \hat{k}_n-\hat{\kappa} \right)^2 \right\}^{1/2}
\]
and
\[
\hat{\rho} (m) =\frac{1}{N-m}
\sum_{n=1}^N \left[\hat{k}_{n+m}-\hat{\kappa})(\hat{k}_n-\hat{\kappa}\right]/\hat{\chi}^2
\]
be the corresponding estimates obtained
from the data $\left\{\hat{k}_n, \;n=1,\ldots, N\right\}$.  We note that $\hat{\kappa}$,
$\hat{\chi}$, and $\hat{\rho} (m)$ have sampling errors proportional to $1/\sqrt{N}$ while $g_1(\theta)$, $g_2(\theta)$, and $g_3(m,\theta)$ have sampling errors proportional to $1/\sqrt{M}$.  The latter sampling errors can be made small by using large values of $M$ in the simulation of Eqs. (\ref{discrete_k_eqn})-(\ref{discrete_p_eqn}).
Setting
\[
J_N^{(2)}(\theta) = \sqrt{\left(g_1(\theta)-\hat{\kappa}\right)^2},\quad J_N^{(3)}(\theta) = \sqrt{\left(g_2(\theta)-\hat{\chi}\right)^2},
\]
and
\[
J_N^{(4)}(\theta) = \sqrt{\frac{1}{L}\sum_{m=1}^{L} \left(g_3(m,\theta)- \hat{\rho} (m)\right)^2},
\]
we estimate $\theta$ by solving
\begin{equation}\label{JObjective}
\min_{\theta} J_N(\theta)=  \min_{\theta} \left\{J_N^{(1)}(\theta)+J_N^{(2)}(\theta)+J_N^{(3)}(\theta)+J_N^{(4)}(\theta)\right\}.
\end{equation}

\begin{table}[H]
\caption{Parameter estimates obtained from quarterly real GDP data (January 1, 1990 -- October 1, 2025).
The critical feedback value at which the deterministic system undergoes a pitchfork bifurcation is
$\gamma_c = 15.12$. The estimated value $\hat{\gamma} = 20.66$ exceeds $\gamma_c$ by approximately 37\%,
placing the system in the regime where the central equilibrium is unstable and two stable outer equilibria exist.
This parameter configuration is consistent with metastable dynamics in $k(t)$ under stochastic perturbations.}
\begin{center}
\begin{tabular}{ccccccccc}
\hline
 & $\hat{\alpha}$ & $\hat{s_1}$ & $\hat{s_2}$ & $\hat{\mu}$ &
$\hat{\gamma}$ & $\hat{\beta}$ & $\hat{a}$ & $\hat{\sigma}$\\
\hline
Estimate & 0.33 & 0.36 & 0.43 & 0.38 & 20.66 & 0.58 & 1.42 & 0.19 \\
\hline
\end{tabular}
\end{center}
\label{DataEstimates}
\end{table}
\section{Monte Carlo Estimates of Stationary Marginal Densities and Business Cycle Time Periods}

\paragraph{Stochastic Metastability and Regime Switching} When $\sigma>0$, the deterministic pitchfork structure described above
persists in a stochastic sense.
For $\gamma>\gamma_c$, the deterministic system possesses two
locally attracting equilibria separated by an unstable saddle at $k=1$.
In the absence of noise, trajectories converge to one of the two
stable equilibria depending on initial conditions.
The introduction of stochastic forcing fundamentally alters this picture.

The diffusion term in Eq. (\ref{discrete_p_eqn}) perturbs the fast variable $p(t)$,
and through the filtered feedback mechanism described in Section~3,
induces fluctuations in $s(t)$ and hence in $k(t)$.
Although these fluctuations are locally small, they accumulate over time
and occasionally drive the system across the unstable separatrix
associated with the saddle equilibrium at $k=1$.
Thus the two deterministic attractors become
\emph{metastable states} of the stochastic system.

From the perspective of stochastic dynamical systems,
the dynamics for $\gamma>\gamma_c$ correspond to a double-well structure
in an effective potential landscape.
Each stable deterministic equilibrium gives rise to a basin of attraction
within which trajectories fluctuate for long periods of time.
Transitions between basins occur through rare noise-induced escape events \cite{LYX}.

This mechanism is classical in the theory of random perturbations of
dynamical systems.
In the small-noise limit, exit times from a basin of attraction
are exponentially distributed to leading order and grow rapidly
as the noise intensity decreases, consistent with
Freidlin--Wentzell large deviation theory \cite{FW}.
Although the present model is not a gradient system,
the geometric mechanism is analogous:
noise drives the system across a dynamically unstable threshold.

\begin{figure}[H]
\begin{center}
\includegraphics[width=12cm]{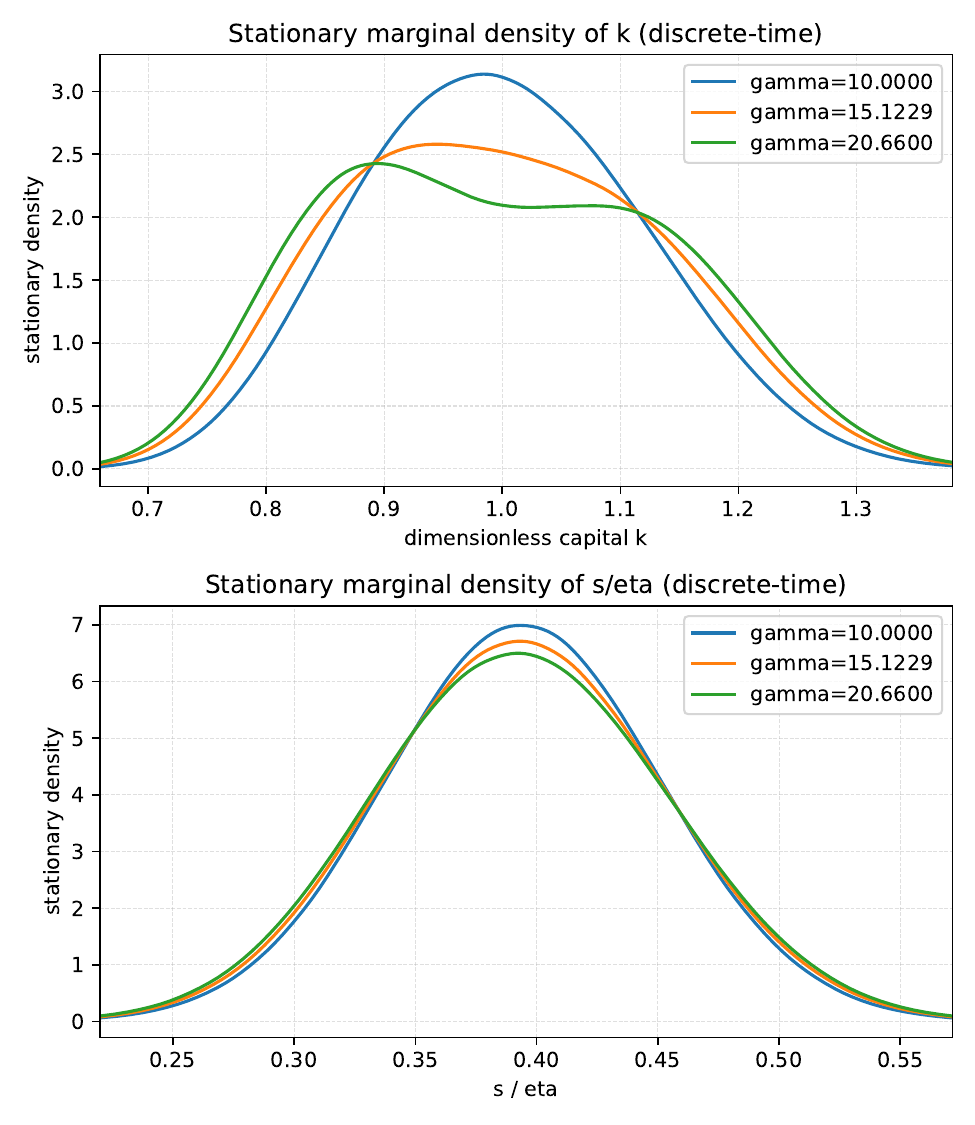}
\caption{Stationary marginal densities of $\left\{k(t), t > 0\right\}$ and
$\left\{s(t), t > 0\right\}$ for values of $\gamma$ below, equal to, and above
$\gamma_c$.}
\label{StationaryDensity}
\end{center}
\end{figure}

\paragraph{Stationary Densities}
Figure \ref{StationaryDensity} shows Monte Carlo estimates of the stationary marginal densities of
$\left\{k(t), t > 0\right\}$ and $\left\{s(t), t > 0\right\}$ of
the system (\ref{discrete_k_eqn})-(\ref{discrete_p_eqn})
using parameter estimates obtained in Section~\ref{DiscreteSection}.
In addition to the value $\gamma= 20.65$, determined from data, and where the
effects of stochastic bifurcation \cite{AT,YW,YLZ,YZD} are clearly evident, we have included the cases $\gamma<\gamma_c$ and $\gamma=\gamma_c$ for comparison.

The metastable structure manifests itself in the stationary marginal
density of $k(t)$.
For $\gamma<\gamma_c$, the density is unimodal and concentrated near $k=1$.
For $\gamma>\gamma_c$, it becomes bimodal, with peaks near the two
stable deterministic equilibria.
The unstable equilibrium at $k=1$ corresponds to the local minimum
between the two modes.

In contrast, the stationary marginal density of $s(t)$ remains unimodal.
This is a direct consequence of the filtering equation
\[
\frac{ds}{dt} = -\beta (s-p),
\]
which smooths high-frequency fluctuations of $p(t)$ and prevents
sharp separation between regimes in the saving rate.
Thus bimodality in $k(t)$ coexists with unimodality in $s(t)$,
reflecting the separation of time scales and the low-pass filtering structure.

\paragraph{Dwell Times and Survival Probabilities}

To quantify regime persistence, we define lower and upper dwell times,
denoted $T_L$ and $T_U$, as the durations for which the trajectory
remains in neighborhoods of the lower ($k<1$) and upper ($k>1$)
deterministic equilibria, respectively.
Operationally, we introduce threshold values $k=0.95$ and $k=1.05$
to separate the two regimes from a central transition region.
A dwell time is defined as the length of a maximal contiguous
subsequence lying entirely on one side of this transition region.

\begin{figure}[H]
\begin{center}
\includegraphics[width=14cm]{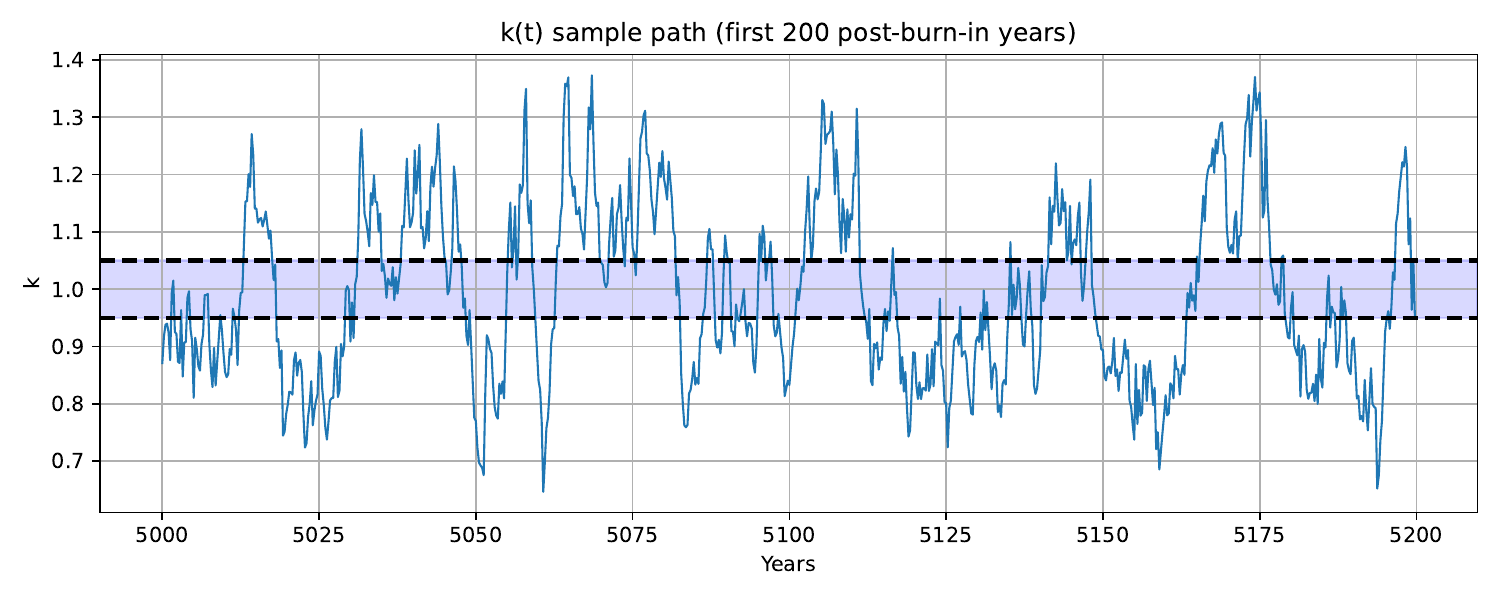}
\caption{Sample path of $k(t)$.  Sample values of $T_L$ and $T_U$ were taken from subsequences of $\left\{k_n, n \ge 0 \right\}$ that lie outside the shaded
transition region.}
\label{Sample_Path}
\end{center}
\end{figure}

\begin{figure}[H]
\begin{center}
\includegraphics[width=14cm]{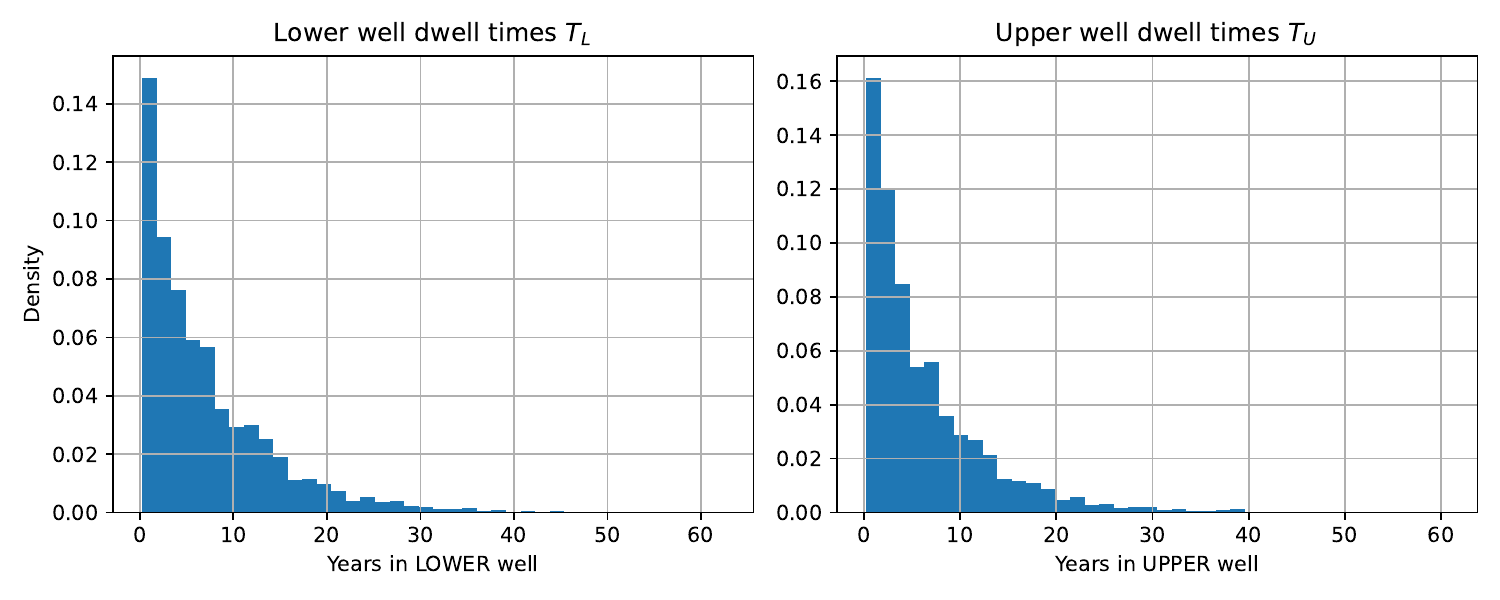}
\caption{Histograms of upper and lower dwell times.}
\label{Dwell_Times}
\end{center}
\end{figure}

After a burn-in period of 5{,}000 years to ensure stationarity,
a 50{,}000-year sample path of the discrete system
(\ref{discrete_k_eqn})--(\ref{discrete_p_eqn}) was generated.
Figure~\ref{Sample_Path} shows a representative 200-year segment of the trajectory.
The shaded region indicates the transition zone separating
the two metastable regimes.
Histograms of the resulting dwell-time samples are shown in Figure~\ref{Dwell_Times}.
Both $T_L$ and $T_U$ exhibit strongly right-skewed distributions.
The median residence times are approximately
5 years for the lower regime and 4.25 years for the upper regime.
Approximately 90\% of dwell times lie below 16.25 years (lower regime)
and 14.75 years (upper regime).

\begin{figure}[H]
\begin{center}
\includegraphics[width=14cm]{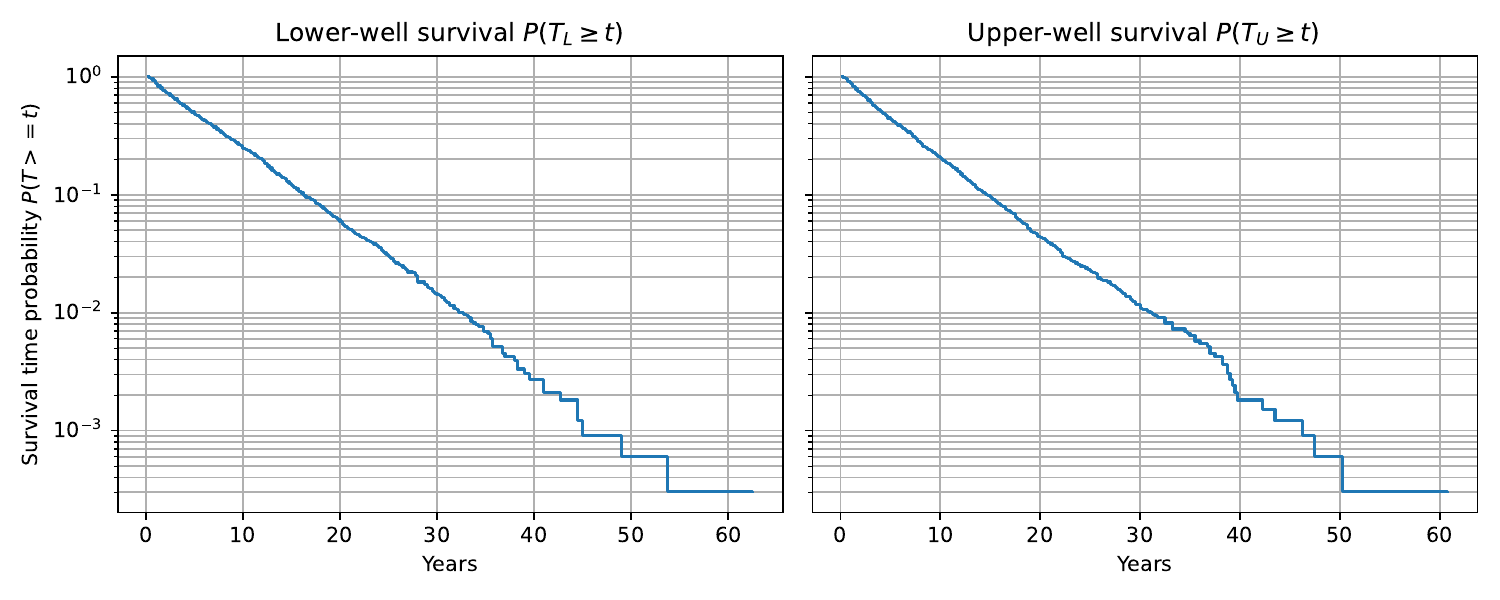}
\caption{Empirical survival functions of lower and upper dwell times}
\label{Survival_Times}
\end{center}
\end{figure}

To examine tail behavior, Figure~\ref{Survival_Times} plots the empirical survival
functions $\mathbb{P}(T_L \ge t)$ and $\mathbb{P}(T_U \ge t)$
on a logarithmic scale.
The near-linear decay in the semi-log plot indicates
approximately exponential tails,
consistent with the theory of noise-induced escape from metastable states.

Within each basin of attraction, trajectories rapidly mix toward
a quasi-stationary distribution.
Escape then occurs as a rare event with approximately constant hazard rate,
leading to survival probabilities of the form
\[
\mathbb{P}(T \ge t) \approx e^{-\lambda t}.
\]
The effective escape rate $\lambda$ depends on both the nonlinear
feedback strength (controlled by $\gamma$) and the noise intensity $\sigma$.
Stronger feedback increases the depth of the effective barrier,
while larger noise intensity lowers it.

\paragraph{Interpretation}

These results confirm the geometric mechanism described earlier.
The deterministic pitchfork bifurcation creates two locally attracting regimes.
Stochastic perturbations transform these into metastable states
and generate random transitions between them.
The observed dwell-time statistics and exponential survival tails
are quantitative signatures of this metastable structure.

Thus business-cycle-like fluctuations in the model arise from
noise-driven transitions between dynamically generated regimes,
rather than from periodic forcing or large exogenous shocks.

\section{Conclusion and Outlook}
\subsection{Summary of Findings}

In this work we have introduced a minimal stochastic growth model in which state-dependent saving behaviour, coupled with bounded multiplicative noise and a natural fast-slow time-scale structure, generates endogenous business-cycle-like dynamics. The model comprises three coupled equations: a Solow-type capital accumulation equation, a linear filtering equation for the saving rate, and a Wright-Fisher-type diffusion for the fast stochastic adjustment process.

The deterministic skeleton of the model exhibits a supercritical pitchfork bifurcation as the feedback gain parameter \(\gamma\) crosses a critical threshold \(\gamma_c\). Below \(\gamma_c\), the balanced-growth equilibrium is unique and stable; above it, two locally attracting equilibria emerge, corresponding to persistent expansion (\(k>1\)) and contraction (\(k<1\)). This bifurcation provides the geometric backbone for regime switching.

When stochastic perturbations are introduced (\(\sigma>0\)), the two deterministic attractors become metastable states. The system spends long periods fluctuating near one regime and occasionally undergoes rare noise-induced transitions to the other. This mechanism reproduces key empirical features of business cycles: persistent expansions and contractions of random duration, bimodal stationary distributions of the detrended capital variable, and right-skewed dwell-time distributions with approximately exponential survival tails.

Using quarterly U.S. real GDP data from 1990 to 2025, we estimated the model parameters via a prediction-error framework augmented by moment-matching constraints. The estimated feedback strength \(\hat{\gamma} = 20.66\) lies well above the critical value \(\gamma_c = 15.12\), placing the system firmly in the bistable regime. Monte Carlo simulations with these parameters yielded stationary densities that are bimodal for \(k(t)\) but unimodal for \(s(t)\) -- a direct consequence of the low?pass filtering of the saving rate. The median dwell times of approximately 5 years for contractions and 4.25 years for expansions align well with the typical durations observed in postwar U.S. business cycles.

\subsection{Outlook and Future Directions}

This work makes three main contributions. First, it demonstrates that realistic business-cycle phenomena can arise endogenously from a low-dimensional nonlinear feedback mechanism, without relying on large exogenous shocks, periodic forcing, or highly disaggregated sectoral structure. Second, it establishes a clear link between deterministic bifurcation theory and stochastic metastability, showing how noise transforms a pitchfork into a double-well potential landscape with observable regime-switching statistics. Third, it provides a practical simulation-based estimation framework that bridges geometric dynamical systems theory with empirical macroeconomic data, making the model amenable to quantitative evaluation. Several extensions and open questions merit further investigation.

\paragraph{Generalisation of the saving function}
The logistic form, while parsimonious, imposes symmetry around \(k=1\). Relaxing this symmetry -- for example, by allowing different saturating levels or asymmetric responsiveness in expansions versus contractions -- could capture more nuanced empirical saving behaviour and potentially yield asymmetric cycle durations, which are often observed in data.

\paragraph{Endogenous noise and volatility clustering}
In the current model, the noise intensity \(\sigma\) is constant. Allowing for state-dependent or time-varying volatility -- e.g., via a stochastic volatility process for the diffusion coefficient -- could generate volatility clustering and time-varying transition rates, further enriching the cycle dynamics.

\paragraph{Multi-country and sectoral extensions}
The present framework is aggregated at the national level. Extending the model to a multi-country setting with trade linkages or to a multi-sector economy could shed light on the propagation of business cycles across regions and industries, and on the role of intersectoral spillovers in regime switching.

\paragraph{Estimation with alternative data and frequencies}
While we used quarterly GDP, the model's time-scale separation suggests that higher-frequency data (e.g., monthly industrial production or saving rates) could provide sharper identification of the fast parameters \(a\) and \(\beta\). Bayesian or likelihood-based estimation methods could also be explored to quantify parameter uncertainty more rigorously.

\paragraph{Large deviations and escape rates}
A more detailed analysis of the noise-induced escape problem -- using Freidlin-Wentzell theory or the action functional -- could yield closed-form approximations for the mean exit times and transition rates as functions of \(\gamma\) and \(\sigma\). This would deepen the theoretical understanding of the regime-switching mechanism and provide testable predictions for the hazard function.

\paragraph{Stochastic bifurcations and tipping points}
The interplay between the deterministic bifurcation and stochastic fluctuations raises interesting questions about tipping phenomena. Investigating how the system responds to gradual changes in \(\gamma\) or \(\sigma\) -- for instance, due to policy shifts or structural changes in saving behaviour -- could inform discussions of critical transitions in macroeconomics and climate-economy interactions.

\paragraph{Incorporating policy and expectations}
The current model treats saving behaviour as a reduced-form function of capital. Embedding forward-looking expectations or policy rules (e.g., interest rate feedback) could make the saving function more microfounded and allow for analysis of monetary and fiscal policy effectiveness in a metastable environment.

In summary, the simple framework developed here offers a fertile ground for exploring endogenous business cycles through the lens of nonlinear dynamics and stochastic processes. Its transparent geometric mechanism, combined with empirical tractability, makes it a valuable benchmark for future theoretical and empirical work on macroeconomic fluctuations. We hope the model stimulates further research into the role of state-dependent behaviour and noise-induced transitions in economic dynamics.

\bigskip

\noindent\textbf{Data Availability}

The numerical algorithms and source code that support the findings of this study are available from the corresponding author upon reasonable request.

\bigskip
\noindent\textbf{Declaration of competing interest}

No author associated with this paper has disclosed any potential or pertinent conflicts which
may be perceived to have impending conflict with this work.

\bigskip
\noindent\textbf{Acknowledgements}

This work is dedicated to Professor James Brannan on the occasion of his 80th birthday.
Shenglan Yuan thanks Pak Hang Chris Lau for insightful discussions on stochastic mechanisms, particularly his theoretical work concerning noise sources and their applications. This work was supported by  Guangdong Basic and Applied Basic Research Foundation (Grant No. 2025A1515012560),  Guangdong Introduction Program (Grant No. 2023QN10X753) and  National Foreign Experts Program (Grant No. 111001819820258003).



\end{document}